\documentclass[11pt,leqno,oneside,letterpaper]{amsart}
\usepackage[letterpaper,left=1.75truein, top=1truein]{geometry}
\usepackage{amsmath,amsfonts,amssymb,amscd,amsthm,amsbsy,epsf,graphics}

\usepackage{color}

\usepackage[colorlinks,linkcolor=blue,citecolor=blue]{hyperref}
\usepackage{epsfig}
\usepackage{graphicx}

\newtheorem{thm}{Theorem}
\newtheorem{cor}[thm]{Corollary}

\newtheorem{prop}[thm]{Proposition}
\newtheorem{conj}[thm]{Conjecture}

\newtheorem{que}[thm]{Question}

\newtheorem{claim}[thm]{Claim}

\theoremstyle{definition}

\def\qed{{\hspace{2mm}{\small $\diamondsuit$}}}

\def\D{{\mathcal D}}

\def\G{{\Gamma}}
 \def\d{{\delta}}

 \def\l{{\lambda}}
 
 \def\m{{\mu}}

   \def\s{{\sigma}}
 
 \def\a{{\alpha}}
 \def\b{{\beta}}
 \def\p{{\partial}}
 \def\r{{\rho}}
 \def\ra{{\rightarrow}}

 \def\lra{{\longrightarrow}}
 \def\g{{\gamma}}
 \def\D{{\Delta}}

 \def\x{{\xi}}
 \def\c{{\mathbb C}}
 \def\z{{\mathbb Z}}
 \def\2{{\mathbb Z_2}}
 
 \def\t{{\tau}}

 \def\da{{\downarrow}}

 \def\sl2{{SL(2,\mathbb C)}}
 
 \def\qed{{\hspace{2mm}{\small $\diamondsuit$}}}

 \def\pf{{\noindent{\bf Proof.\hspace{2mm}}}}

 \def\sk{{{\mbox{\tiny K}}}}

  \def\sl{{{\mbox{\footnotesize  $\mathfrak{L}$}}}}

\begin{document}

\title{Remarks on $SU(2)$-simple knots and $SU(2)$-cyclic $3$-manifolds}

\author{Xingru Zhang}
\address{Department of Mathematics, University at Buffalo, Buffalo, NY 14260}
\email{xinzhang@buffalo.edu}

\maketitle
\vspace{-.6cm}

\begin{center}{\it Dedicated to Steve Boyer on the occasion of his 65th birthday}\end{center}

\begin{abstract}We give some remarks on two closely related issues as stated in the title. In particular we show
that  a  Montesinos knot is  $SU(2)$-simple  if and only if
it is a $2$-bridge knot, extending a result of \cite{Z1} for 3-tangle summand  pretzel knots.
We  conjecture with some evidence that an $SU(2)$-cyclic rational homology $3$-sphere is an $L$-space.
 \end{abstract}

For a knot $K$ in $S^3$, $M_\sk$ will be its exterior and $\m$ a meridian
slope of $K$. Up to a choice of an  orientation for $\m$ and a choice of
the base point for $\pi_1(M_\sk)$, we may also consider $\m$ as an element
 of $\pi_1(M_\sk)$.
  A representation $\r:\pi_1(M_\sk)\ra SU(2)$ is called trace free if
  the trace of $\r(\m)$ is zero (which is obviously well defined).
 An $SU(2)$-representation of $\pi_1(M_\sk)$ is called binary dihedral if its  image is isomorphic to a binary dihedral group.
Note that every binary dihedral representation of $\pi_1(M_\sk)$ is trace free \cite[Proof of Theorem 10]{K}.
  A knot $K$ is called {\it $SU(2)$-simple} if every irreducible trace free
 $SU(2)$-representation of $\pi_1(M_\sk)$ is binary dihedral.

\begin{figure}[!ht]
\centerline{\includegraphics{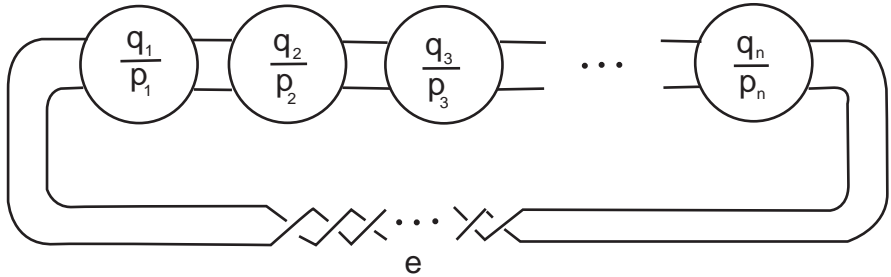}} \caption{A Montesinos link $K(e;q_1/p_1,q_2/p_2,...,q_n/p_n)$.}\label{montesinos}
\end{figure}

 A Montesinos link is usually denoted by $K(e;q_1/p_1,q_2/p_2,...,q_n/p_n)$
where $q_i/p_i$ represents a rational tangle, $|p_i|>1$ and $(q_i,p_i)=1$ for all $i$
(see Figure \ref{montesinos}). By combining the $e$ twists in the figure
with one of the tangles, we may assume that $e=0$, and we will simply write a Montesinos link as
$K(q_1/p_1,\cdots,q_n/p_n)$ and sometimes we refer it as a cyclic tangle sum of $n$ rational tangles.
When $q_i=1, i=1,...,n$, we get a pretzel link.
In \cite{Z1} it was shown that every pretzel knot $K(1/p,1/q,1/r)$, $p,q,r$ pairwise coprime, is not $SU(2)$-simple. In this paper we extend this result to all Montesinos knots of at least three rational tangle summands.

\begin{thm}\label{simple knots}
 Every Montesinos knot $K(q_1/p_1,\cdots,q_n/p_n)$, $n\geq 3$,    is not
 $SU(2)$-simple.
 \end{thm}

For any $SU(2)$-representation $\r:\pi_1(M_\sk)\ra SU(2)$, let $\bar\r:\pi_1(M_\sk)\ra PSU(2)$ be the induced $PSU(2)$-representation.
If $\r:\pi_1(M_\sk)\ra SU(2)$  is a trace free  representation,
then $\r(\m)$ is an order $4$ matrix and $\r(\m^2)=-I$, where $I$ is the identity matrix of $SU(2)$. So $\bar\r:\pi_1(M_\sk)\ra PSU(2)$ factors through the quotient group
$\pi_1(M_\sk)/<\m^2>$, where $<\m^2>$ denote the normal subgroup of $\pi_1(M_\sk)$ generated by $\m^2$.
Let $\Sigma_2(K)$ denote the double branched cover of $(S^3, K)$,
then $\pi_1(\Sigma_2(K))$ is an index two subgroup of
$\pi_1(M_\sk)/<\m^2>$.
It is known that a trace free irreducible representation $\r:\pi_1(M_\sk)\ra SU(2)$ is binary dihedral if and only if  the restriction of $\bar\r$ on $\pi_1(\Sigma_2(K))$ has nontrivial cyclic image \cite[Section I.E]{K}.
For any $2$-bridge knot $K$, $\Sigma_2(K)$ is a lens space and so every irreducible trace free SU(2)-representation of $\pi_1(M_\sk)$ is binary dihedral, that is, every $2$-bridge knot is $SU(2)$-simple.
As a Montesinos knot is a $2$-bridge knot if and only if it has lass than three rational tangle summands,  we have

\begin{cor}A Montesinos knot is $SU(2)$-simple if and only if it is
a $2$-bridge knot.
\end{cor}

By \cite[Corollary 7.17]{KM}, every nontrivial knot in $S^3$ has an irreducible trace free $SU(2)$-representation. It follows that if the double branched cover of a nontrivial knot $K$ is a homology $3$-sphere, i.e. if the knot determinant  $|\D_\sk(-1)|=1$ where $\D_\sk(t)$ is the Alexander polynomial of $K$,  then $K$ is not $SU(2)$-simple.

 A  $3$-manifold $Y$ is called $SU(2)$-cyclic (resp. $PSU(2)$-cyclic) if every $SU(2)$-representation (resp. $PSU(2)$-representation) of $\pi_1(Y)$ has cyclic image. In general $PSU(2)$-cyclic is
 a stronger condition than $SU(2)$-cyclic, that is, $PSU(2)$-cyclic implies $SU(2)$-cyclic but not
 the other way around.
Since $\Sigma_2(K)$ is an $\z_2$-homology $3$-sphere \cite[Cor 3 of 8D]{R}, every $PSU(2)$-representation of $\pi_1(\Sigma_2(K))$ lifts to an $SU(2)$-representation \cite[Page 752]{BZ} and thus
$\Sigma_2(K)$  is $SU(2)$-cyclic if and only if it is $PSU(2)$-cyclic.
So if $\Sigma_2(K)$ is $SU(2)$-cyclic, then $K$ is an $SU(2)$-simple knot.
The following question concerns  the converse.

\begin{que}
Is there  an $SU(2)$-simple knot $K$ in $S^3$ whose double branched cover
$\Sigma_2(K)$ is not $SU(2)$-cyclic (that is, the double branched cover $\Sigma_2(K)$ has irreducible
$PSU(2)$-representations but none of them extend to $M_\sk$)?
\end{que}

One may consider an $SU(2)$-cyclic $3$-manifold  as an $SU(2)$-representation L-space.
The following conjecture suggests that for a rational homology $3$-sphere being $SU(2)$-cyclic  is more restrictive than being a usual
 $L$-space in the  Heegaard Floer homology sense.

 \begin{conj} If   a rational homology $3$-sphere  is
  $SU(2)$-cyclic,  then it is an $L$-space.
\end{conj}

 Here are some evidences for the conjecture.
 Let $K_1=T(p_1,q_1)$ and $K_2=T(p_2,q_2)$ be two torus knots in $S^3$, and let  $M_1$ and $M_2$ be their exteriors. Let $Y(T(p_1,q_1),T(p_2,q_2))$ be the graph manifold obtained by gluing $M_1$ and $M_2$ along their boundary tori
 by an orientation reversing homeomorphism $h:\p M_1\ra \p M_2$
 which identifies the meridian slope in  $\p M_1$ to the Seifert fiber slope in $\p M_2$
 and identifies the Seifert fiber slope in $\p M_1$ with the meridian slope in $\p M_2$.
By \cite[Proposition 5]{Mot} $Y(T(p_1,q_2),T(P_2,q_2))$ has only cyclic $PSL_2(\c)$-representations.
(Although it was assumed in \cite{Mot} that all $p_1, q_1,p_2,q_2$ are positive,  the same argument
 with obvious modification works without this assumption).
Hence $Y(T(p_1,q_1),T(p_2,q_2))$ is $SU(2)$-cyclic.

\begin{prop}\label{L-space}
$Y(T(p_1,q_1),T(p_2,q_2))$ is an $L$-space.
\end{prop}

 \pf We prove  this assertion  by applying
 \cite[Theorem 1.6]{HW}.
By  that theorem, we just need to verify that $h( \mathcal L^{\circ}_{M_1})\cup \mathcal L^{\circ}_{M_2} \cong \mathbb Q\cup \{1/0\}$, where $\mathcal L^0_{M_i}$ is the interior of the set of $L$-space filling slopes of $M_i$, $i=1,2$.
Note that a  general torus knot can be  expressed as $T(p,q)$ with $(p,q)=1$ and $|p|, q\geq 2$.
By \cite[Corollary 1.4]{OS} \begin{equation}\label{L-slopes}{\mathcal L}^0(M_i)=\{\begin{array}{ll}
\mbox{slopes in the open interval  $(p_iq_i-p_i-q_i, \infty)$}, &\mbox{if $p_i>0$,}  \\
\mbox{slopes in the open interval $(-\infty, p_iq_i-p_i+q_i)$}, &\mbox{if $p_i<0$.}\end{array}\end{equation}
Let $\m_i, \l_i$ be the meridian  and longitude  of $K_i$. Note that $p_iq_i$ is the Seifert fiber slope in $\p M_i$.
We have  $h(\m_1)=\m_2^{p_2q_2}\l_2$ and $h(\m_1^{p_1q_1}\l_1)=\m_2$.
Hence for a general slope $m/n$ in $\p M_1$, where $m,n$ are relative prime, $$h(\m_1^m\l_1^n)
=h(\m_1^{m-p_1q_1n}(\m_1^{p_1q_2}\l_1)^n)=(\m_2^{p_2q_2}\l_2)^{m-p_1q_1n}\m_2^n
=\m_2^{p_2q_2(m-p_1q_1n)+n}\l_2^{m-p_1q_1n}.$$
Now suppose $a/b$ is a slope in $\p M_2$, where $a,b$ are relatively prime.
Choose $n=a-p_2q_2b$ and $m=p_1q_1(a-p_2q_2b)+b$, then $m,n$ are relatively prime, $h(m/n)=a/b$,
and \begin{equation}\label{m/n}\frac{m}{n}=p_1q_1+\frac{b}{a-p_2q_2b}=p_1q_1+\frac{1}{\frac{a}{b}-p_2q_2}.\end{equation}

{\bf Case 1}. $p_1>0$ and $p_2>0$.

For any $a/b\notin {\mathcal L}^0(M_2)$, i.e. either  $a/b=1/0$ or $a/b$  is finite and $a/b\leq p_2q_2-p_2-q_2$ by (\ref{L-slopes}),
  choose  correspondingly in $\p M_1$ the slope $m/n=p_1q_1$ or  as in (\ref{m/n})  which
    yields $m/n\geq p_1q_1+\frac{1}{-p_2-q_2}>p_1q_1-1$.
So in either case $m/n\in {\mathcal L}^0(M_1)$ by (\ref{L-slopes}) and $h(m/n)=a/b$, which means $h( \mathcal L^{\circ}_{M_1})\cup \mathcal L^{\circ}_{M_2} \cong \mathbb Q\cup \{1/0\}$ in this case.

{\bf Case 2}. $p_1>0$ and $p_2<0$.

For any $a/b\notin {\mathcal L}^0(M_2)$, we may assume that $a/b$ is finite and so  $a/b\geq p_2q_2-p_2+q_2$
by (\ref{L-slopes}). So $\frac{a}{b}-p_2q_2$ is positive. Choose the slope $m/n$ in $\p M_1$ as in  (\ref{m/n}) which yields  $m/n>p_1q_1$ in this case.
So $m/n\in {\mathcal L}^0(M_1)$ by (\ref{L-slopes}) and $h(m/n)=a/b$.
Thus $h( \mathcal L^{\circ}_{M_1})\cup \mathcal L^{\circ}_{M_2} \cong \mathbb Q\cup \{1/0\}$ holds in this case.

{\bf Case 3}. $p_1<0$ and $p_2>0$.

This case is really Case 2 if we switch $K_1$ and $K_2$.

{\bf Case 4}. $p_1<0$ and $p_2<0$.

For any $a/b\notin {\mathcal L}^0(M_2)$, again we may assume $a/b$ is finite and
so $a/b\geq p_2q_2-p_2+q_2$ by (\ref{L-slopes}). Choose the slope $m/n$ in $\p M_1$ as in (\ref{m/n})
which yields    $m/n\leq p_1q_1+\frac{1}{-p_2+q_2}<p_1q_1+1$.
So $m/n\in {\mathcal L}^0(M_1)$, $h(m/n)=a/b$ and
 we have $h( \mathcal L^{\circ}_{M_1})\cup \mathcal L^{\circ}_{M_2} \cong \mathbb Q\cup \{1/0\}$.

The proof of Proposition \ref{L-space} is now completed.\qed

 It was shown in \cite{Z2} that if $p_1q_1p_2q_2-1$ is odd, then $Y(T(p_1,q_1),T(p_2,q_2))$ is the double branched cover of an alternating knot in
$S^3$, so $Y(T(p_1,q_1),T(p_2,q_2))$  is an $L$-space and the knot in $S^3$ is an $SU(2)$-simple knot
(and is an arborescent knot) but is not a $2$-bridge knot.

There are also examples of hyperbolic rational homology $3$-spheres which are $SU(2)$-cyclic \cite{C}.
These examples are also double branched covers of alternating knots in $S^3$ and thus are $L$-spaces.
These alternating knots are thus $SU(2)$-simple but are not arborescent.

{\bf Proof of Theorem \ref{simple knots}}.
Let $K=K(q_1/p_1,\cdots,q_n/p_n)$ be a Montesinos knot with $n\geq 3$. We need to show that
$\pi_1(M_\sk)$ has an irreducible trace free $SU(2)$-representation which is not binary
dihedral. Here is an outline of how the proof goes.
We show  that the double branched cover $\Sigma_2(K)$ has an irreducible
$PSU(2)$-representation $\bar\r_0$ which can be extended to an $PSU(2)$-representation $\bar\r$ of $\pi_1(M_\sk)$ up to conjugation. This  $PSU(2)$-representation $\bar\r$ lifts to an $SU(2)$-representation
$\r$ of $\pi_1(M_\sk)$ which is  automatically trace free. Since $\bar\r_0$ is an irreducible representation,
$\r$ is not  binary dihedral. The existence of $\bar\r_0$ is provided by \cite{B}.
We first apply some ideas from \cite{Mat} to show that $\bar\r_0$
extends to a unique $PSL_2(\c)$-representation $\bar\r$ of $\pi_1(M_\sk)$.
Then we further show that this $\bar\r$ is conjugate to  an $PSU(2)$-representation by applying
some results from \cite{HP}\cite{CD}.

Now we give the details of the proof.
For a finitely generated group $\G$, $\bar R(\G)=Hom(\G, PSL_2(\c))$
denotes the $PSL_2(\c)$ representation variety of $\G$ and
$\bar X(\G)$ the $PSL_2(\c)$ character variety of $\G$.
Let $t:\bar R(\G)\ra \bar X(\G)$ be the map which sends a representation $\bar\r$ to its character
$\chi_{\bar\r}$.
We shall write an element in $PSL_2(\c)$ as $\bar A$ which is the image of an element $A$ in $SL_2(\c)$
under the quotient map $SL_2(\c)\ra PSL_2(\c)$ and for convenience we sometimes call elements
in $PSL_2(\c)$ as matrices.
For any $\bar A\in PSL_2(\c)$ define $tr^2(\bar A)=(trace(A))^2$ which is obviously well defined.
Recall that the character $\chi_{\bar \r}$ of an $PSL_2(\c)$-representation $\bar\r$ is the function
$\chi_{\bar\r}:\G\ra \c$ defined by $\chi_{\bar\r}(\g)=tr^2(\bar \r(\g))$.

A character $\chi_{\bar\r}$ is real if $\chi_{\bar \r}(\g)\in \mathbb R$
for all $\g\in\G$.
If we consider $\bar X(\G)$  as an algebraic subset in $\c^n$ (for some $n$),
then real characters of $\bar X(\G)$   correspond to real points of $\bar X(\G)$, i.e.
points of  $\bar X(\G)\cap \mathbb R^n$.
If $\s:\c^n\ra\c^n$ (for each $n\geq 1$)  denotes the operation of coordinatewise  taking complex conjugation,
 then any complex affine algebraic set $Y$  in $\c^n$ defined over $\mathbb Q$ is invariant under $\s$
and the set of real points of $Y$ is precisely the fixed point set of $\s$ in $Y$.
Note that $\bar R(\G)$ and  $\bar X(\G)$ are both algebraic sets defined over $\mathbb Q$ and that the map $t:\bar R(\G)\ra \bar X(\G)$ is an algebraic map defined over $\mathbb Q$, we thus have the following
commutative diagram of maps:
$$\begin{array}{ccc}
\bar R(\G)&\stackrel{\s}{\lra}&\bar R(\G)\\
\da t&&\da t\\
\bar X(\G)&\stackrel{\s}{\lra}&\bar X(\G).
\end{array}$$
It follows that
$\s(\chi_{\bar\r})=\chi_{\s(\bar\r)}$.

A representation $\bar \r\in \bar R(\G)$ is irreducible if the image of $\bar \r$
cannot be conjugated into the set \{$\bar A$; $A$ upper triangular\}.
Two irreducible representations in $\bar R(\G)$ are conjugate
if and only if  they have the same character.

If $W$ is a compact manifold, $\bar R(W)$ and $\bar X(W)$ denote $\bar R(\pi_1 W)$ and $\bar X(\pi_1W)$
respectively.

\begin{figure}[!ht]
\centerline{\includegraphics{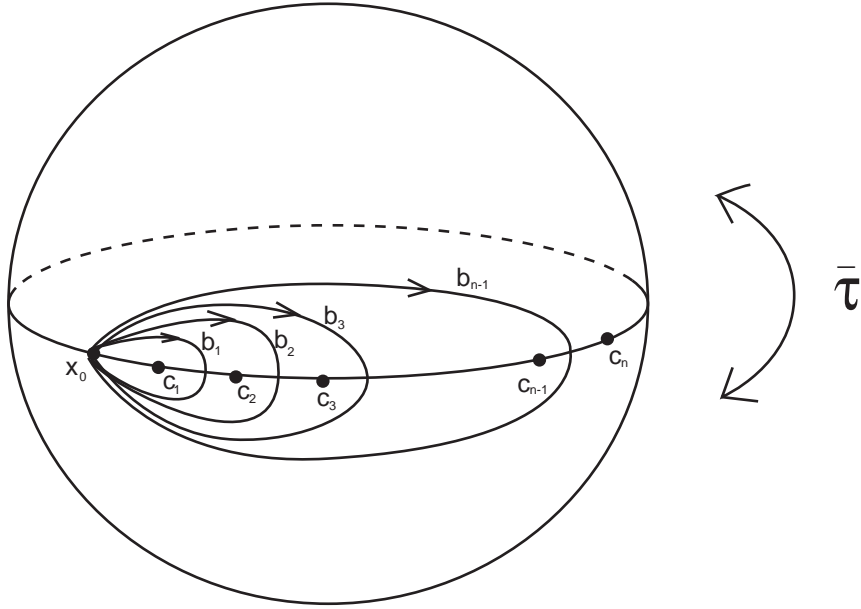}} \caption{The orbifold $S^2(p_1,\cdots, p_n)$, its involution $\bar\t$ and the generating set $b_1,...,b_{n-1}$ for $\D(p_1,...,p_n)$, where $x_0$ is the base point and $c_1,..., c_n$ are cone points of orders $p_1,...,p_n$ respectively. }\label{orbifold}
\end{figure}

Let $K=K(q_1/p_1,\cdots,q_n/p_n)$ and $M=M_\sk$. We may assume that all $p_i$ are positive
(by changing the sign of $q_i$ if necessary).
Let $p:\tilde M\ra M$ be the  $2$-fold cyclic  covering
and let $\tilde\m=p^{-1}(\m)$ which is a connected simple closed essential curve in
$\p \tilde M$ which double covers $\m$.
Then $p_*:\pi_1(\tilde M)\ra \pi_1(M)$ is an injection
and we may consider $\pi_1(\tilde M)$ as an index two normal
subgroup of $\pi_1(M)$, in which $\tilde \m=\m^2$.
Dehn filling $\tilde M(\tilde \m)$ of $\tilde M$ with the slope $\tilde \m$ is the double branched cover $\Sigma_2(K)$ of $(S^3, K)$.
The covering involution $\t$ on $\tilde M$ extends to one
on $\tilde M(\tilde \m)$ which we still denote by $\t$.
Montesinos proved in \cite{Mo1}\cite{Mo2}
 that $\tilde M(\tilde\m)$  admits a Seifert fibering
 invariant under the covering involution $\t$, the base orbifold of the Seifert fibred space is $S^2(p_1,...,p_n)$
which is the $2$-sphere with $n$ cone points of orders $p_1,...,p_n$, and  $\t$ descends
down to an involution $\bar \t$ on $S^2(p_1,...,p_n)$ which is a reflection in a circle passing through all the cone points (see Figure \ref{orbifold}).

 We denote the orbifold fundamental group of $S^2(p_1,...,p_n)$ by $\D(p_1,...,p_n)$
 which has the following presentation:
 $$\D(p_1,...,p_n)=<a_1,...,a_{n}; a_i^{p_i}=1, i=1,...,n, a_1a_2\cdots a_n=1>.$$
It was shown in \cite[Section 3.3]{Mat} that when $n=3$ any irreducible
$PSL_2(\c)$-representation of $\pi_1(\tilde M)$ which factors
through $\pi_1(\tilde M(\tilde \m))$ has a unique extension
to $\pi_1(M)$. Note that this extended representation can be lifted to a trace free $SL_2(\c)$-representation of $\pi_1(M)$.
We shall slightly extend this result to the following

\begin{prop}\label{extension}Let $\delta$ be the composition of the
three  quotient homomorphisms
$$\pi_1(\widetilde M)\ra \pi_1(\widetilde M(\tilde\m))\ra
\D(p_1,p_2,...,p_n)\ra \D(p_1,p_2,p_3).$$
Let $\phi:\D(p_1,p_2,p_3)\ra PSL_2(\c)$ be any irreducible
representation.
Then $\bar\r_0=\phi\circ \d$ has a unique extension to $\pi_1(M_\sk)$.
\end{prop}

\pf The proof for uniqueness is verbatim as that given in \cite{Mat} on page 38-39.
We need to note that  as $K(q_1/p_1,\cdots,q_n/p_n)$  is a knot
at most one of $p_i$'s is even. So \cite[Lemma 2.4.9]{Mat}
still applies to our current case.

\begin{claim}\label{A} There will be an extension $\bar\r$ if and only if there is
$\bar A\in PSL_2(\c)$ such that ${\bar A}^2=\bar I$ (where $I$ is the identity matrix of $SL_2(\c)$) and $\bar A\bar\r_0(\b)\bar A^{-1}=\bar\r_0(\m\b\m^{-1})$
for all $\b\in \pi_1(\tilde M)$.\end{claim}

Again this claim can be proved verbatim as that of \cite[Claim 3.3.2]{Mat}.

So to finish the proof of Proposition \ref{extension}, we just need to find an $\bar A\in PSL_2(\c)$ with the properties
stated  in Claim \ref{A}, which is what we are going to do in the rest of the proof of Proposition \ref{extension}.
Recall that $\tilde M(\tilde\m)$ is the Dehn filling of $\tilde M$ with a solid torus
$N$ whose meridian slope is identified with the slope $\tilde\m$.
The core circle of $N$ is the fixed point set of $\t$ in $\tilde M(\tilde\m)$.
Let $D$ be a meridian disk of $N$ such that the fixed point of $\t$ in $D$ (the center point of $D$)
is disjoint
from the singular fibers of the Seifert fibred space $\tilde M(\tilde\m)$.
Choose a point  $\tilde x$ in $\p D$ and let $\tilde x_0$ be the center point of $D$.
Then  arguing  as on \cite[Page 40]{Mat} we have  the following commutative
diagram:
$$\begin{array}{ccccccc}
\pi_1(\tilde M, \tilde x)&\lra&\pi_1(\tilde M(\tilde \m), \tilde x)&\lra
\pi_1(\tilde M(\tilde \m),\tilde x_0)&\lra&\D(p_1,....,p_n)\\
&&&&&&\\
\da (\cdot)^\m&&&\da \bar \t_*&&\da \bar\t_*\\
&&&&&&\\
\pi_1(\tilde M,\tilde x)&\lra&\pi_1(\tilde M(\tilde \m),\tilde x)&\lra\pi_1(\tilde M(\tilde \m),\tilde x_0)&\lra &\D(p_1,....,p_n)
\end{array}$$
where $(\cdot)^\m:\pi_1(\tilde M,\tilde x)\ra\pi_1(\tilde M,\tilde x)$ corresponds to the conjugation action by $\m$, i.e.  $(\beta)^\m=\m\beta\m^{-1}$ and $\D(p_1,...,p_n)$ is the orbifold fundamental group of $S^2(p_1,...,p_n)$ whose base point is the image $x_0$ of the point $\tilde x_0$
under the quotient  map $\tilde M(\tilde \m)\ra S^2(p_1,...,p_n)$.

Figure \ref{orbifold} shows the generating set $b_1, ...,b_{n-1}$ of the orbifold fundamental group $\D(p_1,...,p_n)$ of $S^2(p_1,...,p_n)$.
In fact we have
$$a_1=b_1, a_2=b_1^{-1}b_2, a_3=b_2^{-1}b_3, \cdots, a_{n-1}=b^{-1}_{n-2}b_{n-1}, a_n=b_{n-1}^{-1}$$
and conversely
$$b_1=a_1, b_2=a_1a_2, b_3=a_1a_2a_3, \cdots, b_{n-1}=a_1a_2\cdots a_{n-1}, b_{n-1}=a_n^{-1}.$$
Obviously from Figure \ref{orbifold} the induced isomorphism $\bar \t_*:\D(p_1,...,p_n)\ra \D(p_1,...,p_n)$
sends $b_i$ to $b_i^{-1}$, $i=1,...,n-1$.
So we have
$\bar\t_*(a_1)=a_1^{-1}$, $\bar\t_*(a_2)=b_1b_2^{-1}=a_1a_2^{-1}a_1^{-1}$,
$\bar\t_*(a_3)=b_2b_3^{-1}=a_1a_2a_3^{-1}a_2^{-1}a_1^{-1}$, $\cdots$,
$\bar\t_*(a_{n-1})=b_{n-2}b_{n-1}^{-1}=a_1a_2\cdots a_{n-2}a_{n-1}^{-1}a_{n-2}^{-1}\cdots a_2^{-1}a_1^{-1}$,
 $\bar\t_*(a_n)=\bar\t_*(b_{n-1}^{-1})=b_{n-1}=a_n^{-1}$.
Since the quotient homomorphism
 $$\begin{array}{l}\D(p_1,\cdots, p_n)=<a_1,..., a_n; a_i^{p_i}=1, i=1,...,n, a_1a_2\cdots a_n=1>
\\\lra \D(p_1, p_2, p_3)=<\bar a_1, \bar a_2, \bar a_3; \bar a_i^{p_i}=1, i=1,2,3, \bar a_1\bar a_2\bar a_3=1>\end{array}$$
sends $a_i$ to $\bar a_i$, $i=1,2,3$, and send $a_i$ to $1$, $i=4,...,n$.
we see that $\bar\t_*$ descents to an isomorphism
$\bar\t_{\#}:\D(p_1,p_2, p_3)\ra \D(p_1,p_2, p_3)$ such that
$\bar\t_{\#}(\bar a_1)=\bar a_1^{-1}$, $\bar\t_{\#}(\bar a_2)=\bar a_1\bar a_2^{-1}\bar a_1^{-1}$,
$\bar\t_{\#}(\bar a_3)=\bar a_1\bar a_2\bar a_3^{-1}\bar a_2^{-1}\bar a_1^{-1}$ and we have the following commutative diagram:
$$\begin{array}{ccccccccc}
\pi_1(\tilde M)&\lra&\pi_1(\tilde M(\tilde \m))&\lra &\D(p_1,....,p_n)&\lra&\D(p_1,p_2,p_3)&\stackrel{\phi}{\lra} PSL_2(\c)\\
&&&&&&&&\\
\da (\cdot)^\m&&\da \t_*&&\da \bar \t_*&&\da \bar\t_{\#}&&\\
&&&&&&&&\\
\pi_1(\tilde M)&\lra&\pi_1(\tilde M(\tilde \m))&\lra &\D(p_1,....,p_n)&\lra&\D(p_1,p_2,p_3)&\stackrel{\phi}{\lra} PSL_2(\c)
\end{array}$$
So $\bar\t_\#(\bar a_1\bar a_2)=\bar a_2^{-1}\bar a_1^{-1}=(\bar a_1\bar a_2)^{-1}$.
Since $\D(p_1,p_2,p_3)$ is generated by $\bar a_1, \bar \a_2$, we see by applying  \cite[Lemma 3.1]{BZ}
that $\phi$ and $\phi\circ\bar \t_\#$ have the same $PSL_2(\c)$ character.
(In fact if $\phi(\bar a_1)=\bar A_1$ and $\phi(\bar a_2)=\bar A_2$, then $\phi(\bar a_1\bar a_2)=\bar A_1\bar A_2=\overline{A_1A_2}$, $(\phi\circ\bar\t_\#)(\bar a_1)
 =(\bar A_1)^{-1}=\overline{A_1^{-1}}$, $(\phi\circ\bar\t_\#)(\bar a_2)=\bar A_1(\bar A_2)^{-1}(\bar A_1)^{-1}=\overline{A_1 A_2^{-1}A_1^{-1}}$ and $(\phi\circ\bar\t_\#)(\bar a_1\bar a_2)=(\bar A_2)^{-1}(\bar A_1)^{-1}=\overline{(A_1A_2)^{-1}}$.
 Now let $F_2$ be the free group on two generators $\xi_1$ and $\xi_2$.
 Let $\r_1$ and $\r_2$ be the $SL_2(\c)$ representations of $F_2$ defined by $\r_1(\xi_i)=A_i$, $i=1,2$,
 and $\r_2(\xi_1)=A_1^{-1}$, $\r_2(\xi_2)=A_1A_2^{-1}A_1^{-1}$. Then one can easily verify that $tr(\r_1(\xi_1))=tr(\r_2(\x_1))$,
 $tr(\r_1(\xi_2))=tr(\r_2(\x_2))$ and $tr(\r_1(\xi_1\xi_2))=tr(\r_2(\x_1\xi_2))$. So
 \cite[Lemma 3.1]{BZ} applies.)
So $\phi$ and $\phi\circ\bar \t_\#$ are conjugate $PSL_2(\c)$ representations, that is, there is $\bar A\in PSL_2(\c)$ with $\bar A\phi\bar A^{-1}=\phi\circ\bar\t_\#$.
Combining this with the definition of $\bar\r_0$ and the last commutative diagram, we see that
$\bar A\bar\r_0(\b)\bar A^{-1}=\bar\r_0(\m\b\m^{-1})$ for each $\b\in \pi_1(\tilde M)$.
The proof of Proposition \ref{extension} is now finished.\qed

Now by \cite{B}, every triangle group $\D(p_1,p_2,p_3)$  has an irreducible
$SO(3)\cong PSU(2)$-representation.
Therefore there is an irreducible representation $\bar\r_0$ as given in Proposition \ref{extension}
with its image contained in $PSU(2)$.
So the character $\chi_{\bar\r_0}$ of $\bar\r_0$ is real valued.
Let $\bar\r$ be the unique extension of $\bar\r_0$ to $\pi_1(M)$ as guaranteed by Proposition \ref{extension}.

\begin{claim}The character $\chi_{\bar\r}$ of $\bar\r$ is real valued.
\end{claim}

Suppose otherwise. Recall that  $\s: \bar X(M)\ra \bar X(M)$ is the operation of taking complex conjugation and  a character is real valued if and only if it is a fixed point of $\s$.
So $\chi_{\bar\r}\ne \s(\chi_{\bar\r})=\chi_{\s(\bar\r)}$ are two different characters of irreducible
representations and
thus $\bar\r$ and $\s(\bar \r)$ are non-conjugate representations.
But $\chi_{\bar\r_0}=\s(\chi_{\bar\r_0})=\chi_{\s(\bar\r_0)}$ and $\bar\r_0$ is irreducible.
Hence $\bar\r_0$ and $\s(\bar\r_0)$ are conjugate representations, that is, there is $\bar B\in PSL_2(\c)$ such that $\bar\r_0=\bar B\s(\bar\r_0)\bar B^{-1}$.
Hence $\bar\r$ and $\bar B\s(\bar\r)\bar B^{-1}$ are non-conjugate representations which have the
same restriction on $\pi_1(\tilde M)$.
We get a contradiction with  Proposition \ref{extension}.

By \cite[Lemma 10.1]{HP} an $PSL_2(\c)$-character $\chi_{\bar\r}$ of a finitely generated group is real valued if and only if
the image of $\bar\r$ can be conjugated into $PSU(2)$ or $PGL_2( \mathbb R)$.
So  our current representation $\bar\r$ can be conjugated into $PSU(2)$ or $PGL_2(\mathbb R)$.
If it can be conjugated into $PSU(2)$, then we are done because this conjugated representation lifts
to a trace free $SU(2)$-representation which is not binary dihedral.
So suppose that $\bar \r$ is conjugate to an $PGL_2(\mathbb R)$-representation $\bar\r'$.
As noted  in \cite{HP} right after Lemma 10.1, $PGL_2(\mathbb R)\subset PGL_2(\c)\cong PSL_2(\c)$ has two components, the identity component is $PSL_2(\mathbb R)$
and the other component consists of matrices of determinant  $-1$
(which in $PSL_2(\c)$ are represented by matrices with entries in $\c$
with zero real part). Considering the action of $PSL_2(\c)$ on hyperbolic space
$\mathbb H^3$ by orientation preserving isometries, the group $PGL_2(\mathbb R)$ is the stabilizer of a total geodesic plane
in $\mathbb H^3$, and an element of $PGL_2(\mathbb R)$ is contained in $PSL_2(\mathbb R)$ if and only if it preserves the orientation of the plane.
Since $\pi_1(\tilde M)$ is the unique index two normal subgroup of $\pi_1(M)$, $\pi_1(\tilde M)$ is generated by elements $\g^2$, $\g\in \pi_1(M)$.
Therefore the image of the restriction $\bar\r_0'$ of $\bar\r'$ on $\pi_1(\tilde M)$ consists of elements preserving the orientation
of the total geodesic  plane mentioned above, i.e. the image of $\bar\r_0'$ is contained in $PSL_2(\mathbb R)$. So the image of $\bar \r_0$  can  be conjugated into $PSL_2(\mathbb R)$. But the image of $\bar\r_0$ is contained in $PSU(2)$.
\cite[Lemma 2.10]{CD} says that if an $PSU(2)$-presentation can be conjugated into an $PSL_2(\mathbb R)$-representation, then it is a reducible representation. So our $\bar\r_0$ is a reducible representation. We arrive at a contradiction, which completes the proof of Theorem \ref{simple knots}.

\def\bysame{$\underline{\hskip.5truein}$}

\end{document}